\documentclass[12pt,oneside,english,shortlabels]{amsart}

\usepackage[T1]{fontenc}
\usepackage[latin9]{inputenc}
\usepackage{color}
\usepackage{babel}
\usepackage{amstext}
\usepackage{amssymb}
\usepackage{esint}
\usepackage[all]{xy}
\usepackage[unicode=true,pdfusetitle,
 bookmarks=true,bookmarksnumbered=false,bookmarksopen=false,
 breaklinks=false,pdfborder={0 0 0},backref=false,colorlinks=true]
 {hyperref}
\hypersetup{
 linkcolor=blue}
\usepackage{breakurl}

\makeatletter

\usepackage{amsthm}
\usepackage{enumitem}		

\usepackage{mathrsfs} 
\providecommand{\noopsort}[1]{}

\let\cal\mathcal

\def\11{{\mathbf 1}}

\def\CC{{\mathbb C}}

\def\NN{{\mathbb N}}

\def\QQ{{\mathbb Q}}
\def\RR{{\mathbb R}}

\def\ZZ{{\mathbb Z}}

\def\l{\lambda}

\theoremstyle{plain}

\newtheorem{thm}[subsubsection]{Theorem}
\newtheorem{lem}[subsubsection]{Lemma}
\newtheorem{cor}[subsubsection]{Corollary}
\newtheorem{prop}[subsubsection]{Proposition}

\theoremstyle{definition}
\newtheorem{defn}[subsubsection]{Definition}

\theoremstyle{remark}
\newtheorem{remark}[subsubsection]{Remark}

\theoremstyle{remark}
\newtheorem{exampl}[subsubsection]{Example}

\theoremstyle{remark}
\newtheorem{notation}[subsubsection]{Notation}

\def\bee{\begin{exampl}}
\def\eee{\end{exampl}}
\def\bn{\begin{notation}}
\def\en{\end{notation}}
\def\br{\begin{remark}}
\def\er{\end{remark}}
\def\bp{\begin{prop}}
\def\ep{\end{prop}}
\def\bpr{\begin{proof}}
\def\epr{\end{proof}}
\def\bt{\begin{thm}}
\def\et{\end{thm}}
\def\be{\begin{equation}}
\def\ee{\end{equation}}
\def\bl{\begin{lem}}
\def\el{\end{lem}}
\def\bc{\begin{cor}}
\def\ec{\end{cor}}
\def\bd{\begin{defn}}
\def\ed{\end{defn}}

\numberwithin{equation}{subsection}

\author{Georges Comte}
\thanks{}
\address{Univ. Grenoble Alpes, Univ. Savoie Mont Blanc, CNRS, LAMA, 73000 Chamb\'ery, France}
\email{georges.comte@univ-smb.fr}
\urladdr{http://gc83.perso.sfr.fr/}

\author{Yosef Yomdin}
\thanks{This research was supported by the ISF, Grant No. 779/13, by FSMP and by ANR DEFIGEO.
The authors would like to thank University of Savoie Mont Blanc and the Weizmann Institute
 for the hospitality during part of the
research for this paper. The authors are very grateful to the anonymous referee for its careful reading and for
many insightful suggestions.}
\address{Department of Mathematics, The Weizmann Institute of Science,
Rehovot 76100, Israel}
\email{yosef.yomdin@weizmann.ac.il}
\urladdr{http://www.wisdom.weizmann.ac.il/~yomdin/}

\title[Zeroes and  rational points of analytic functions]
{Zeroes and  rational points of analytic functions}

\begin{document}

\dedicatory{\today}
\maketitle

\begin{abstract}
For an analytic function $f(z)=\sum_{k=0}^\infty a_kz^k$ on a neighbourhood of a closed disc $D\subset  \CC$, we give assumptions, in terms of the Taylor coefficients $a_k$ of $f$, under which the number of intersection points  of the graph $\Gamma_f$ of $f_{\vert D}$ and algebraic curves of degree $d$ is polynomially bounded in $d$. In particular, we show these assumptions are satisfied for random power series, for some explicit classes of lacunary series, and for solutions of algebraic differential equations with coefficients
and initial conditions in $\QQ$. As a consequence, for any function $f$ in these families,  $\Gamma_f$ has less than $\beta \log^\alpha T$ rational points of height at most $T$, for some $\alpha, \beta >0$.

\end{abstract}

\tableofcontents

\section{Introduction}\label{Sec:Introduction }
\setcounter{equation}{0}

Let $D_R\subset \CC$  be a closed disc, centred at the origin and of radius $R$, let $f:D_R\to \CC$ be an analytic function on a neighbourhood of $D_R$, and for any $d\ge 1$, let us denote by ${\cal P}_d$ the subspace
of polynomials of $\CC[X,Y]$ of degree at most $d$.  Assuming that for any polynomial $P\in \CC[X,Y] $, $P(z,f(z))$ is not identically zero, or in other words, that $f$ is transcendental, the quantity
$$
Z_d(f)=\underset{P\in {\cal P}_d}{\sup}\#\{z\in D_R; P(z,f(z))=0\},
$$
(where zeroes are counted with multiplicity)
is finite.
The integer $Z_d(f)$ is also the maximum number of intersection points between the graph $\Gamma_f$ of $f$ and algebraic curves of degree at most $d$. Thus any bound of $Z_d(f)$ will be called in the sequel a B\'ezout bound for $f$, since such a
bound is called a B\'ezout bound in the case of proper intersection of two algebraic curves. Moreover in the algebraic case a polynomial bound in the degrees involved
holds.

\smallskip

Bounding $Z_d(f)$ is very closely related to the study of ``Doubling'', or ``Bernstein - type'' inequalities for restrictions of polynomials $P$ as above to the graph $\Gamma_f$ of $f$ (see, for instance, \cite{Com.Pol,Roy.Yom}). Through this connection, many methods and results in Potential Theory, on one side, and in Approximation Theory, on the other side, become relevant, starting with the basic results of \cite{Sad1,Sad2}. Still, we try to keep the references to the minimum.

\smallskip 

This article studies the growth of $Z_d(f)$, as $d$ goes to infinity; more specifically, we provide certain assumptions on the transcendental analytic function $f$, under which $f$ has, like in the algebraic case, a B\'ezout bound polynomial in $d$.

The asymptotic behaviour of $Z_d(f)$ as $d$ goes to infinity was studied for different classes of functions, with different tools.
 In fact,  $Z_d(f)<\infty $ generally holds when $f$ is definable in an o-minimal structure expanding the real field	
\footnote{Note that, on the opposite side, by \cite{GwKuPa}, a polynomial B\'ezout bound for $f$ does not imply that $f$ is polynomial or even definable in some o-minimal structure.}.
But in this very general situation the behaviour of $f$ with respect to algebraic curves of growing degree is difficult to predict. In the analytic case, considered in the present paper, the following observation is instructive: for any $\zeta\in ]0,1[$ there is analytic function  $f$ such that for a sequence of degrees $d$ going to infinity, $Z_d(f)\ge e^{d^\zeta}$ (see, for instance, \cite[Example 7.5]{Pi2004}, \cite{Su2002} or \cite{Su2006}, and inequality (\ref{eq: lot of zeroes}) of Remark \ref{rmk: Pila Wilkie is asymptotically sharp} below).
On the other hand, one also knows that for any analytic function $f$, $Z_d(f)$ is bounded from above by a polynomial in $d$ of degree $2$, for a certain sequence of degrees $d$ going to infinity (see \cite[Theorem 1.1]{Com.Pol2007}) and the asymptotic of this upper bound is best possible (see for instance \cite[Corollary 2.6]{Com.Pol})


A polynomial in $d$ B\'ezout bound for $f(z)=e^z$ was obtained in \cite{Tij}. For entire functions $f$ of positive order (under some additional conditions) such bounds were obtained in \cite{Com.Pol,Com.Pol2007}
\footnote{  For the first time instances of analytic functions (with lacunary Taylor series at the origin) having polynomial B\'ezout bound with prescribed growth were provided
in \cite[Corollary 6.2]{Com.Pol}.}
and in \cite{BoJo15b, Bru2008}. For the Riemann zeta function, and for some other specific functions accurate polynomial bounds were obtained in \cite{Bos.Bru.Lev,Ma11}. Very recently, for solutions of certain types of algebraic differential equations, a polynomial in $d$ B\'ezout bound was obtained in \cite{BinarXiv}.


Our approach to the problem of bounding $Z_d(f),$ for $f(z)=\sum_{k=0}^\infty a_kz^k,$ is based on a detailed algebraic study of the Taylor coefficients $a_k$ of $f$. It follows a long research line, starting with a classical work of Bautin (\cite{Bau1,Bau2}). Bautin's discovery was that for analytic families $f_\lambda(z)=\sum_{k=0}^\infty v_k(\lambda)z^k$, with $\lambda$
-- a parameter ranging in a finite dimensional space, and $v_k(\lambda)$ -- polynomials in $\lambda$, the number of zeroes can be bounded in terms of the polynomial ideals $I_r=\{v_0(\lambda)
,\ldots,v_r(\lambda)\}$
(the Bautin ideals).
 This approach was further developed in many publications (see, for instance, \cite{Bru2003,Bru2015,Fri.Yom, Roy.Yom,Yom1} and references therein). In particular, the case of linear families (that is to say $v_k(\lambda)$ is linear in $\lambda$) where computations are reduced to linear algebra, was considered in \cite{Yom}. In the present paper we concentrate on linear families of the form $P(z,f(z))$. Notice that some notions introduced below play important role also in Hermite-Pad\'e approximations (see, for instance, \cite{Bak.Lub,Nik.Sor}, and references therein, and Remark \ref{remark:transc.seq} below). We plan to present some instances of this important connection separately.


We will consider polynomials $P(z,y)=\sum_{0\le i,j \le d} \lambda_{i,j}z^iy^j$.
In this case the parameter vector $\lambda = (\lambda_{i,j}$)
ranges in the space $\CC^m$, with $ m=(d+1)^2$ and we
have
$$
f_\lambda(z)=P(z,f(z))=\sum_{k=0}^\infty v_k(\lambda)z^k,
$$
with $v_k(\lambda)$ linear forms in the parameter $\lambda$. The specific coefficients of $v_k(\lambda)$ can be explicitly written through the Taylor coefficients $a_k$ of $f$ (see Section \ref{Sec: Bezout bounds for transcendental analytic curves} below). For each $k$ we can write
$$
\begin{pmatrix}
v_0(\lambda) \\
\vdots \\
v_k(\lambda)
\end{pmatrix}
= M_k
 \begin{pmatrix}
\lambda_1 \\
\vdots \\
\lambda_m
\end{pmatrix},
$$
with a
$(k+1)\times m$
 matrix $M_k$. We define the Bautin index $b_d=b_d(f)$ as the minimal $k$ for which the rank of $M_k$ is equal to $m$, and call $M_{b_d}=M_{b_d}(f)$ the $d$-th Bautin matrix of $f$
(see Definition \ref{de: Bautin matrix}).
The Bautin index $b_d$ is also the maximal, with respect to $\lambda$, multiplicity at $0$ of $f_\lambda(z)$
(see Proposition \ref{prop:max multiplicity and Bautin index}).
  Finally, we denote by $\delta_d=\delta_d(f)>0$ the maximum of the absolute value of all non-zero minor determinants of size $m\times m$ of the Bautin matrix $M_{b_d}$. In this way we associate to each transcendental analytic function $f(z)=\sum_{k=0}^\infty a_kz^k,$ two sequences $(b_d)_{d\ge 1}$, and $(\delta_d)_{d\ge 1}$. For hypertranscendental $f$ (that is when all the derivatives $f^{(r)}$, $r\ge 0$, are algebraically independent) we also define a sequence of multiplicities $(\eta_d)_{d\ge 1}$
%
in the following way
(see also Definition \ref{def:Hyper.Transc.Ind}).

 Consider the Taylor expansion of $f(z)=\sum_{k=0}^\infty a_k(u)(z-u)^k$ at points $u$ near the origin. Using these expansions we construct, as above, the square matrices $M_m(u), \ m=(d+1)^2$. Assume that the determinants $\Delta_d(u)$ of $M_m(u)$ do not vanish identically in $u$, and let $\eta_d(f)$ be the multiplicity of zero of $\Delta_d(u)$ at $u=0$. We then have
$$
\Delta_d(u) = \alpha_d u^{\eta_d}+O(u^{\eta_d+1}), \hbox{ with } \alpha_d\ne 0.
$$

The three sequences $(\delta_d), (b_d), (\eta_d)$
are defined algebraically, through a finite amount (depending on $d$) of the Taylor coefficients $a_k$ of $f$, and they present natural transcendence measures of $f$. In the present paper we investigate some basic properties of these sequences, and their role in bounding $Z_d(f).$ All the  three sequences can be defined for formal power series $f(z)$ with   coefficients in an arbitrary field, although here we work only with complex (or real) analytic functions on $D_R$.

\smallskip

\noindent The following are our main results:


\begin{itemize}

\item[(1)] We bound $Z_d(f)$ in terms of $b_d(f)$ and $\delta_d(f)>0$
(Proposition \ref{Prop:Est.c}).
 In particular, for $b_d(f)\leq R(d),$ and $\delta_d(f)\geq e^{-S(d)},$ with $R,S$ polynomials in $d$, we have $Z_d(f)\leq T(d),$ with $T$ also a polynomial in $d$.

\item[(2)] If the Taylor coefficients $a_k$ of $f$ are rational, we define $h_l(f)$ as the maximal denominator of $a_k$, $k=0,1,\ldots,l$. We bound from below $\delta_d$ in terms of $b_d$ and $h_{b_d}$
(Proposition  \ref{prop:rat.coef.delta}).
In particular, for $b_d(f)\leq R(d),$ and $h_l(f)\leq e^{S(l)},$ with $R$ and $S$ polynomials, we show that $\delta_d(f)\geq e^{-U(d)},$ and hence $Z_d(f)\leq T(d),$ with $T,U$ also polynomials in $d$
(Theorem~\ref{thm: Bound through the transcendency index sequence}).
As an example, we show that this is the case for solution of algebraic ODE's with rational coefficients and rational initial values
(Theorem~\ref{thm:ODE}).
This gives another proof of one of results in \cite{BinarXiv}.

\item[(3)]
 If the Taylor coefficients $a_k$ at $0$ of a hypertranscendental function $f$  are rational and satisfy $h_l(f)\leq e^{S(l)},$ with  $S$ polynomial, and if there exists
 a polynomial $R$ such that
 $\eta_d\le R(d)$, $d\ge 1$, then we show that $Z_d(f)\leq T(d),$ with $T$ also polynomial.
(Theorem~\ref{thm: Bezout Bound through polynomial hypertranscendence}).

\item[(4)] We consider a class of lacunary series $f$, similar to the one considered in \cite{ComPol2003}. We show that for this class the Bautin index $b_d(f)$ can be explicitly estimated, as well as the Bautin matrices and the determinants $\delta_d(f).$ On this base we obtain examples of both polynomial and non-polynomial growth of $Z_d(f)$
(Theorem~\ref{thm: Bezout Bound for lacunary series}).

\item[(5)] Clearly, for a series $f$ with random Taylor coefficients $a_k,$ the square Bautin matrices $M_m(u), \ m=(d+1)^2$ are non-degenerate. Hence $b_d(f)=m=(d+1)^2.$ We show that with probability one the determinants $\Delta_d$ of these matrices satisfy $\Delta_d\geq e^{-U(d)}$
(Theorem~\ref{thm: Bound for Delta at random f }),
and therefore $Z_d(f)\leq T(d),$ with $T,U$ polynomials in $d$
(Corollary~\ref{cor: Bezout at random}).

\item[(6)]
Analytic functions on compact domains with polynomial B\'e\-zout bounds have the following remarkable Diophantine property: the number of  rational points of height $\le T$ in their graph is bounded by
a power of $\log T$. It means that they have few rational points of given height, since a sharp upper bound for this number, for analytic functions, is in $C_\epsilon T^\epsilon$, for any $\epsilon>0$. Therefore our assumptions also provide families of transcendental sets with few rational points, as explained in Section \ref{Sec:Rat.Points}
(see Theorem \ref{thm: bound for rational points of bounded height }).
\end{itemize}
\smallskip

The paper is organized as follows. In section \ref{Sec:Lin.fam} we consider general linear families $f_\lambda$, we introduce the definitions used throughout the paper, and we give a general bound for $Z_d(f)$ in terms of the Bautin index and the determinant $\delta$. In section \ref{Sec: Bezout bounds for transcendental analytic curves} we prove the results (1) to (4). In section \ref{Sec: Bautin determinant for random series} we prove (5). Finally in Section \ref{Sec:Rat.Points} we give direct Diophantine applications of our analytic B\'ezout theorems, as mentioned in $(6)$.



%

\section{Linear families of analytic functions}\label{Sec:Lin.fam}
\setcounter{equation}{0}


The beginning of this Section is based on \cite{Yom}.

Let us denote by $D_R$ the disc $\{ z\in  \CC, \vert z\vert \le R \}$,
let $\psi:(D_R,0)\to(\CC^n,0)$ be an analytic curve, and let

\be\label{eq:curve}
\psi(t)=\sum^{\infty}_{k=1}a_k z^k, \ a_k\in \CC^n,
\ee
be the Taylor expansion of $\psi$ at $0\in\CC$. We assume that the series (\ref{eq:curve}) converges in a neighbourhood of the disk $D_R=\{ z\in  \CC, \vert z\vert \le R \}$, and that
$||\psi(t)||\leq A$ for any $z\in D_R$.

For $\lambda=(\lambda_1,\ldots,\lambda_m)\in\CC^m$, let
$$ Q_{\lambda}=\sum^m_{i=1}\lambda_iQ_i,$$
with $Q_j:( \Omega ,0) \to\CC$, $j=1,\ldots,m$, analytic functions in a neighbourhood  $\Omega$ of the polydisc of
$\CC^n$ of radius $A$, and   bounded there by $B$.

In what follows we are interested in linear families of analytic functions of the form
\be\label{eq:f.lambda}
f_{\lambda}(z)=Q_{\lambda}(\psi(z))=\sum^m_{i=1}\lambda_i Q_i(\psi(z)).
\ee
Write
\be\label{eq:curve2}
Q_i(\psi(z))=\sum^{\infty}_{k=0}c^i_k z^k.
\ee
By our assumptions $Q_i(\psi(z))\leq B$ on $D_R$, $i=1,\ldots,m$.
Hence, by Cauchy's estimates, we have
\be\label{eq:curve2.1}
|c_k^i|\leq \frac{B}{R^k}.
\ee
For the Taylor development of the linear family $f_\lambda$ we have
\be\label{eq:curve1}
f_{\lambda}(z)=\sum^{\infty}_{k=0}v_k(\lambda) z^k, \ \ v_k(\lambda)=\sum^m_{i=1}c^i_k\lambda_i .
\ee
Thus the coefficients $v_k(\lambda)$ of the series (\ref{eq:curve1}) are linear forms in $\lambda$ and by (\ref{eq:curve2.1}), for any $k$
\be\label{eq:curve3}
\vert v_k(\lambda)\vert\leq \frac{m B |\lambda|}{R^k} ,
\ee
where $|\lambda|=\max \{|\lambda_1|,\ldots, |\lambda_m| \}$.

Now we associate to the family $f_\lambda$ an integer $b=b(f_\lambda)$. Let, for $i\in \NN$,  $L_i\subseteq\CC^m$ be the linear
subspace of $\CC^m$ defined by the equations $v_0(\lambda) = \ldots = v_i(\lambda)=0$. We have
$$
L_0\supseteq L_1\supseteq \ldots L_i\supseteq \ldots .
$$
Hence on a certain step $b$ this sequence stabilizes
$$
L_{b-1}\supsetneqq L_b=L_{b+1}=\cdots =L.
$$
Note that $\lambda\in L_b \Longleftrightarrow v_k(\lambda)=0,$
for $k=0,1, \ldots$
\bd
\label{de: Bautin index}
We call this number $b$ the {\sl Bautin index of the family $f_\lambda$} (see
\cite{Bau1,Bau2}).
\ed
\br\label{rem:transcendeal functions and Bautin index}
Of course the dimension of the subspaces $L_i\subset \CC^m$, $i\in \NN$, is at most
$m$ and if $z\mapsto f_\lambda(z)$ is not identically zero for all $\lambda$, $\dim(L_b)\le m-1$. More accurately, if we assume that for $\lambda \ne 0$ the function $f_\lambda(t)$ does not vanish identically, then we have $L_b=\{0\}$.
 But the dimension of $L_i$ can drop at most by one at each step, thus necessarily in this situation $b\ge m-1$.
 This will be the case if $f(z)$ is a transcendental function for the family of polynomials of degree at most $d$, that is for $\psi(z)=(z,f(z))$ and $Q_{i,j}=X^iY^j$ with $i,j\le d$, a classical case mainly considered in the following sections.

On the other hand the Bautin index may be as big as wished. In general one cannot explicitly find the Bautin index of $f_\lambda$, since the moments, when the dimension of $L_i$ drops, are usually difficult to determine.
\er
A first characterization of the Bautin index of the family $f_\lambda$ is the following.
\bp\label{prop:max multiplicity and Bautin index}
Let us assume that for $\lambda\not=0$ the function
$z\mapsto f_\lambda(z)$ is not identically zero, and let
us denote by $\mu$ the maximal multiplicity, with respect to the parameters $\lambda\not=0$,
 of the Taylor series at the origin of $f_\lambda(z)$. Then $\mu=b$.
\ep

\bpr
   There exists a  parameter $\lambda$ such that $f_\lambda(z)$
has multiplicity $\mu$, that is such that
$$f_\lambda(z)=v_\mu(\lambda) z^\mu +v_{\mu+1}(\lambda)z^{\mu+1}+\cdots, $$
with $v_\mu(\lambda)\not=0$.
Therefore $\lambda\in L_{\mu-1} \setminus L_\mu$. It follows that $L_\mu \subsetneqq L_{\mu-1}$, and  $b\ge \mu$. On the other hand,  since no parameter $\lambda\not=0$ can cancel $v_0, \ldots, v_\mu$ in the same time,
$L_\mu=\{0\}=L_{\mu+1}=\cdots$, and thus $b\le \mu$.
\epr

\br
The system $v_0(\lambda)=\cdots = v_{m-2}(\lambda)=0$, with $m$ parameters, always
having a non-zero solution, one sees that the maximal multiplicity of $f_\lambda$ is at
least $m-1$. Proposition \ref{prop:max multiplicity and Bautin index} then implies that $b\ge m-1$,
that was already noted in Remark \ref{rem:transcendeal functions and Bautin index}.
\er

\br\label{rem:asymptotic number of zeroes and maximal multiplicity}
Under the assumption that for $\lambda\not=0$ the function
$z\mapsto f_\lambda(z)$ is not identically zero, and counting zeroes with multiplicity, from
Proposition \ref{prop:max multiplicity and Bautin index} one observes that
$$ \lim_{r\to 0}\underset{\lambda\not=0}{\max}\#\{z\in D_r, f_\lambda(z)=0\} \ge b. $$
We give hereafter in Theorem \ref{thm:zroes.baut} more accurate relations between the Bautin index and the number of
zeroes of the family $f_\lambda$, showing in particular that the above inequality is an equality.
\er

Since $L$ is defined by $v_0(\lambda)=\cdots =v_b(\lambda)=0$, any linear form $\ell(\lambda)$, which vanishes on $L$, can be expressed as a linear
combination of $v_0,\ldots,v_b$. Now a basis $(v_{i_1}, \ldots, v_{i_\sigma})$, $i_1, \ldots, i_\sigma \in \{0,\ldots, b\}$, of the space of linear forms vanishing on $L$ being chosen among the elements of the family
$(v_0,\ldots, v_b)$, there exists a constant $\tilde c >0$, depending on this basis, such that for any $\ell$ with
$\ell_{\vert L} \equiv 0$, we have

\be\label{eq:curve4}
\ell(\lambda)=\sum^\sigma_{j=1}\mu_jv_{i_j}(\lambda), \ \mu_j\in\CC, \ and \ \vert\mu_j\vert\leq \tilde c\Vert\ell\Vert, j=1,\ldots,\sigma.
\ee
where for $\ell(\lambda)=\sum^m_{i=1}\alpha_i\lambda_i$, $\Vert\ell \Vert =\max_i\vert\alpha_i\vert$.

\bn\label{no:c}
We denote by $c=c(f_\lambda)>0$ the minimum of the constants $\tilde c$ satisfying (\ref{eq:curve4}).

\en
An effective estimation of $c$ is difficult, in general. However, if the
Bautin index $b(f_\lambda)$ is known, $c(f_\lambda)$ can be estimated via a finite computation in terms of the Taylor coefficients of $\psi$
and of $Q_i$ (see Proposition \ref{Prop:Est.c} below).

Now let a family $f_\lambda$ be given, and let the Bautin index $b(f_\lambda)$ and the constant $c(f_\lambda)$ of this family be
defined as above.

\bt\label{thm:zroes.baut}
Zeroes being counted with multiplicity, we have the following uniform bounds with respect to the parameter $\lambda$, when $f_{\lambda}\not\equiv 0$.
\begin{enumerate}
\item The maximal number of zeroes of  $f_{\lambda}$ in the disk $D_{\frac{R}{4}}$ is at most
$$ 5 b \log \big(4+2c(b+1)\frac{B}{R^b} \big) \ \hbox{ if }  R\le 1,$$
$$ \hbox{ and } \
 5 b \log \big(4+2c(b+1)B \big) \ \hbox{ if }  R\ge 1,$$

\item This maximal number of zeroes of $f_\lambda$ is at most $b$
in $D_\rho$, where
$$
\rho=\frac{R}{e^{10b+2}\max(2,c(b+1)B\max(\frac{1}{R},1)^b) }.
$$

\end{enumerate}
\et

\bpr
 For any $\lambda\in\CC^m$, and for any $j\geq b+1$, by  the definition of $c$ and from the bound (\ref{eq:curve2.1}) (see Theorem 1.1 of \cite{Yom} and the last inequality in its proof), we have
\be\label{eq:T.dom}
|v_j(\lambda)|R^j \leq c(b+1)B\max(\frac{1}{R},1)^b\ \underset{i=0,\ldots,b}{\max} |v_i(\lambda)|R^i.
\ee
 Then, the bounds to prove on the number of zeroes
  are consequences of (\ref{eq:T.dom}) and
\cite[Lemma 2.2.1, Theorem 2.1.3]{Roy.Yom}.
\epr
\br\label{rem:asymptotic number of zeroes equals maximal multiplicity}
Note that under the assumption that for $\lambda\not=0$ the function
$z\mapsto f_\lambda(z)$ is not identically zero, and counting zeroes with multiplicity, from
Theorem  \ref{thm:zroes.baut} and Remark \ref{rem:asymptotic number of zeroes and maximal multiplicity} one has
$$ \lim_{r\to 0}\underset{\lambda\not=0}{\max}\#\{z\in D_r, f_\lambda(z)=0\} = b. $$
 This infinitesimal maximal number of zeroes is
called the {\sl cyclicity of $Q_\l$ on $\psi$} (see \cite{Yom}).
\er
In what follows we develop explicit bounds on the number of zeroes of $f_\lambda$ in $D_{\frac{R}{4}}$ in terms of the coefficients $c^i_k$, with
$i=1,\ldots,m, \ k=0,\ldots,b$, i.e., ultimately, in terms of the Taylor coefficients of $\psi$ and of $Q_i$ up to the order $b$.
By Theorem \ref{thm:zroes.baut}, this amounts bounding
the constant $c(f_\l)$ introduced in Notation \ref{no:c}. Note also that such a bound is given by a bound on the coefficients $\mu_j, \ldots, \mu_\sigma$ of the system (\ref{eq:curve4}), since in this system one can consider only linear forms $\ell$ of norm $1$.
%
Let us denote the dimension of the stabilized subspace $L=L_b$ by $s=m-\sigma \ (\leq m-1)$.
All the information we need, as we will see, is encoded in  the rank $\sigma$ matrix
$ M=(c^i_k), \ k=0,\ldots,b, \ i=1,\dots,m$, with $b+1$ lines and $m$ columns.
With our notation, $M=M_b$ is defined by
$$  \begin{pmatrix}
v_0(\lambda) \\
\vdots \\
v_b(\lambda)
\end{pmatrix}
= M_b
 \begin{pmatrix}
\lambda_1 \\
\vdots \\
\lambda_m
\end{pmatrix}.
$$
\bd
\label{de: Bautin matrix}
With the above notation, the $(b+1)\times m$ matrix $M=M_b$
is called the {\sl Bautin matrix of the family $f_\lambda$}.
The matrix $M_b$ is the matrix of the linear map
$\Lambda $ sending elements of the vector space $\cal Q$ spanned by
the analytic functions $Q_1,\ldots, Q_m$ (assumed to be linearly independent) to the space of
$b$-jets at the origin of analytic functions. The map $\Lambda$ is the composition of the linear maps $\widetilde{\Lambda}:{\cal Q}\to
\CC\{z\}$ sending $Q_\lambda$ to $f_\lambda$ with the linear map
$j^b_0:\CC\{z\}\to \CC^{b+1}$ of $b$-jets at the origin, $b$ being the first order of jets at the origin such that $\dim ( j^b_0(\widetilde{\Lambda}({\cal Q}))=\dim (\widetilde{\Lambda}({\cal Q})).$
\ed
\bn\label{not: small delta}
We will denote by $\delta>0$ the maximum of the absolute value of all  non-zero minor
determinants of size $\sigma\times \sigma$ of the Bautin matrix $M$.
\en


\bp\label{Prop:Est.c}
Let $f_\lambda$ be given as above. Then
$$
c(f_\lambda) \le
 \sigma\frac{ (B\sqrt \sigma)^{\sigma-1} }{\delta R^{\beta(\sigma-1)} },
 \  \hbox{ where } \beta=b \  \hbox{ if }  R\le 1 \
 \hbox{ and  } \  \beta=\frac{\sigma}{2} \  \hbox{ if }  R\ge 1. $$
In turn, in the disk $D_{\frac{R}{4}}$, the maximal number $Z(f_\lambda)$ of zeroes of the family $f_\lambda$, with respect to the parameter $\lambda$,  satisfies

$$
Z(f_\lambda) \le 5b
\log (
4+ 2(b+1) \sqrt\sigma
\frac{ (B\sqrt\sigma)^{\sigma} } { \delta R^{ b \sigma }}
)
\ \hbox{ if } R\le 1,
$$
$$
\hbox{ and } \ Z(f_\lambda) \le 5b
\log (
4+ 2(b+1) \sqrt\sigma
\frac{ (B\sqrt\sigma)^{\sigma} }{ \delta R^{ \frac{\sigma}{2}(\sigma-1) } }
)
\ \hbox{ if } R\ge 1,
$$

\ep
\bpr
Let $\hat{M}$ be a submatrix of $M$ of size $\sigma\times\sigma$, with the absolute value of the determinant equal to $\delta$, according to notation \ref{not: small delta}. $\hat{M}$ is
obtained from $M$ by the choice of the lines $i_1 < \cdots < i_\sigma$
of $M$ (corresponding to the choice $(v_{i_1},\ldots, v_{i_\sigma})$ for a basis of the space of linear forms cancelling on $L$), and the choice of certain $\sigma$ columns, say,
the first $\sigma$ columns in $M$.
When
$\Vert \ell \Vert \le 1$,  (\ref{eq:curve4})
gives the linear system
$$  \begin{pmatrix}
\alpha_1 \\
\vdots \\
\alpha_\sigma
\end{pmatrix}
= {}^t\hat{M}
 \begin{pmatrix}
\mu_1 \\
\vdots \\
\mu_\sigma
\end{pmatrix}
$$
with $\vert \alpha_j\vert \le 1$, $j=1, \ldots, \sigma$.
Therefore, by the Cramer rule, each $\mu_j$ satisfies
\be\label{eq:curve6}
|\mu_j|\leq \frac{ \sigma \hat{\delta}}{\delta},
\ee
where $ \hat{\delta}$ is the maximum of the absolute values of $(\sigma-1)\times (\sigma-1)$ sub-minors of $\hat{M}$. Next, by (\ref{eq:curve2.1}), we have
$$
|c_{i_j}^i|\leq \frac{B}{R^{i_j}}\leq \frac{B}{R^b} \  \hbox{ if } R\le 1$$
$$\hbox{ and } \
|c_{i_j}^i|\leq \frac{B}{R^{i_j}}\leq \frac{B}{R^j} \ \hbox{ if } R\ge 1.
$$
Consequently, the length of the $j$-th row-vectors in $(\sigma-1)\times (\sigma-1)$ sub-minors of $\hat{M}$ does not exceed
$$
\frac{B\sqrt{\sigma-1}}{R^b} \  \hbox{ if } R\le 1,
\ \hbox{ and } \
\frac{B\sqrt{\sigma-1}}{R^j} \  \hbox{ if } R\ge 1
$$
Interpreting the determinant as the volume of the span of its
row-vectors,
we conclude that

$$
\sigma\hat{\delta}
\leq
\sigma \frac{ (B\sqrt{\sigma-1})^{\sigma-1} }{R^{b(\sigma-1)} }
\le
 \sigma\frac{ (B\sqrt \sigma)^{\sigma-1} }{R^{b(\sigma-1)} } \  \hbox{ if } R\le 1
$$
$$
\hbox{ and } \ \sigma\hat{\delta}
\leq
\sigma \frac{ (B\sqrt{\sigma-1})^{\sigma-1} }{R^{\frac{\sigma}{2}(\sigma-1)} }
\le
 \sigma\frac{ (B\sqrt \sigma)^{\sigma-1} }{R^{\frac{\sigma}{2}(\sigma-1)} } \  \hbox{ if } R\ge 1.
$$
This bound, combined with (\ref{eq:curve6}) and Theorem \ref{thm:zroes.baut}, completes the proof of Proposition \ref{Prop:Est.c}.
\epr

\br\label{rmk: centred at 0}
From the beginning of this section, we have assumed that the analytic map $\psi$ is defined on a disc centred at the origin and that $\psi(0)=0$. This choice is harmless since, in case $\psi$ is defined on a ball centred at $a$, one can consider
$\phi(z)= \psi(z+a)-\psi(a)$ and the bounds given in Theorem
\ref{thm:zroes.baut} and Proposition \ref{Prop:Est.c} for $\phi$ and the family $Q_i(w+\psi(a))$
are the same bounds for $\psi$ and the family $Q_i(w)$ when Taylor series are considered at $a$ instead of $0$.
\er
\br\label{rmk: bounded by 1 and defined on D4}
A classical application of bounds given in Proposition \ref{Prop:Est.c} is for analytic plane curves $\psi(z)=(z,f(z))$ and
the family $Q_i(X,Y)$ of two variables monomials of total degree at most some integer $d$.
In this case we consider that $f$ is given by its Taylor series
at the origin and that this series converges on the disc $D_R$ of radius $R$ and centred at the origin, by Remark \ref{rmk: centred at 0}.
The curve $\psi(z)$ is the standard parametrization of the graph of $f$. In this setting we will
provide bounds for the number of zeroes of $P_d(z,f(z))$   on $D_R$ that are uniform with respect to coefficients of polynomials $P_d$ of degree at most $d$.
Now note that in this situation one can only consider functions $f$ that are bounded by $1$ on $D_R$, since a uniform bound on the number of zeroes of $P_d(z,\frac{1}{N} f(z))$, where $N$ bounds $f(z)$ on $D_R$, provides a uniform bound  on the number of zeroes of $P_d(z,f(z))$.
In the same way one can as well assume for simplicity that $f$ is analytic on the unit disc, up to applying the bounds provided
by Proposition \ref{Prop:Est.c} to  the new function $g(z)=f(Rz)$, since a uniform bound on the zeroes of
$P_d(w,g(w))$ on $D_1$ is a uniform bound on the zeroes of
$P_d(z,f(z))$ on $D_R$.
Of course the same rescaling effects apply in the same way
for the family $X^iY^j, i,j\le d$, of two variables monomials. Nevertheless those reductions are not always possible for any family of
analytic functions $Q_j$ and for any analytic curve $\psi$.

Finally when one restricts to analytic functions bounded
by $1$ on $D_1$ and for the family $Q_{i,j}=X^iY^j, i+j\le d$
 or  for the family $X^iY^j, i,j\le d$
 of two variables monomials, one can take $1$ for $B$ as bound for $Q_{i,j}(z,f(z))$ on $D_1$.
Note that in this case, from Proposition \ref{Prop:Est.c} we deduce the following statement.
\er
\bc\label{cor : Bezout bound in case B=R=1 }
 Let $f_\lambda$ be given as above with $R=1$ and $B=1$, then in the disk $D_{\frac{1}{4}}$ the maximal number $Z(f_\lambda)$ of zeroes of the family $f_\lambda$, with respect to the parameters $\lambda$,  satisfies
$$
Z(f_\lambda) \le 5b \log (4+ 2(b+1)\frac{ e^{\sigma\log\sigma}}{\delta}   ) .
$$
\ec

\bpr
After taking $B=1$ and $R=1$ in the last inequality of Proposition \ref{Prop:Est.c},
observe that $\sqrt{\sigma}^{\sigma +1}=e^{(\sigma +1)\frac{1}{2}
\log \sigma}\le e^{\sigma\log\sigma}$.
\epr

\section{B\'ezout bounds for transcendental analytic curves}\label{Sec: Bezout bounds for transcendental analytic curves}
\setcounter{equation}{0}

\subsection{Families of polynomials}
From now on we consider the classical  case of the family of polynomials with degree at most a fixed integer $d$. That is to say, for a given analytic function
$f:(D_1,0)\to (\CC,0)$

we want to bound the number of zeroes
of $P_d(z,f(z))$, for $P_d$ a two variables polynomial of degree at most $d$.
Considering Remarks \ref{rmk: centred at 0} and \ref{rmk: bounded by 1 and defined on D4}, for simplicity we assume that $f$ sends $0$ to $0$, is given by its converging Taylor series on the unit disc $D_1$, and is bounded by $1$ on $D_1$.

In fact, we consider the problem in the following (essentially equivalent) form
\footnote{ \ The family $ \sum_{0\le i,j\le d} \lambda_{j,i}z^if^j(z)$ considered here presents the advantage to have a slightly more symmetric Bautin matrix then the family $ \sum_{i+j\le d} \lambda_{j,i}z^i f^j(z)$. Of course one can easily deduce
the Bautin matrix of the second family from the Bautin matrix of the first family, and the problem could be considered only for the family
$ \sum_{i+j\le d}\lambda_{j,i} z^i f^j(z)$ as well. In particular note that
$Z_d(f)\le \cal Z_d(f)\le Z_{2d}(f)$, and therefore a  polynomial bound for the maximal number
of zeroes of one family gives rise to a   polynomial bound for the maximal number
of zeroes of the other family.   }:
give a bound for
$$\cal Z_d(f):=\underset{p_j\in \CC[z], \deg p_j\le d}{\max}  \#   \{ z\in D, \sum_{j=0}^d p_j(z) f^j(z)=0\}$$
for some disc $D\subset D_1$ centred at the origin.
Denoting
$ p_j(z)= \sum_{i=0}^d\lambda_{j,i}z^i,$
the sum
$\sum_{j=0}^d p_j(z) f^j(z)$
has the form
$ \sum_{0\le i,j\le d} \lambda_{j,i}z^if^j(z).$

With the notation of Section \ref{Sec:Lin.fam}, we have to consider
the family of monomials $Q_{i,j}=X^iY^j$, $i,j\le d$, and the analytic function $\psi(z)=(z,f(z))$, therefore here the number
of parameters is $m=(d+1)^2$.
Let us denote by ${\cal Q}_d$ the vector space spanned by the monomials $Q_{i,j}=X^iY^j$, $i,j\le d$.
In order to guarantee $\cal Z_d(f)<\infty$ for any $d$ as soon as some $p_j$ is not zero,
we assume that $f$ is transcendental, that is no non-zero polynomial restricted to the graph of $f$ vanishes identically.
By Remark \ref{rem:transcendeal functions and Bautin index}, the Bautin index $b=b_d$ of this family satisfies
$b\ge m-1=d^2+2d$ and the corresponding Bautin matrix $M=M_b$ has at least $m=(d+1)^2$ lines (and exactly $m$ columns). Furthermore the dimension $s$ of the space $L_b$ is $0$ and therefore the dimension
$\sigma$ of the space generated by $v_0, \ldots, v_b$ is $m$.
The matrix $M$ is the matrix of the linear map
$\Lambda : {\cal Q_d}\to \CC^{b+1}$ (in the basis of monomials) sending an element  $P(X,Y)\in {\cal Q}_d$ to the $b$-jet
$j^b_0(P(z,f(z))$ at the origin of the analytic map $P(z,f(z))$, and, with the notation of Definition \ref{de: Bautin matrix}, $b$ is the first index of
jets such that $\dim(j^b_0(\widetilde{\Lambda}({\cal Q}_d)))=\dim(\widetilde{\Lambda}({\cal Q}_d))=(d+1)^2$. 

\bn\label{no: coefficients a_j^i}
For $i,j\in \NN$, we denote by $a_i^j$ the $i$th Taylor coefficient at the origin of the $j$th power $f^j$ of $f$. Namely,
$ a_i^j=\frac{1}{i!}(f^j)^{(i)}(0)$.
\en

With this notation, a direct computation shows that

$$
\hskip-3mm M=\begin{pmatrix}
1 & 0 & \rule[0.5ex]{6mm}{0.1mm} &  0 & a_0^1 & 0 &\rule[0.5ex]{6mm}{0.1mm} & 0 &   \cdots & a_0^d & 0 & \rule[0.5ex]{6mm}{0.1mm} & 0 \\
 \rule[0mm]{0.1mm}{10mm} &  &  &    \rule[0mm]{0.1mm}{10mm}&  \rule[0mm]{0.1mm}{10mm} & & &   \rule[0mm]{0.1mm}{10mm} & &   \rule[0mm]{0.1mm}{10mm} & & &       \rule[0mm]{0.1mm}{10mm} \\
0& &  & 0 & a_{d-1}^1 &  a_{d-2}^1 & & 0 & \cdots & a^d_{d-1} & &  &0 \\
0   & &  & 1  &  a_{d}^1 & a_{d-1}^1 & & a_0^1 &  \cdots & a^d_{d}   & &  & a^d_0       \\
0   & &  & 0 & a_{d+1}^1 & a_d^1 &  & a_1^1 &  \cdots & a^d_{d+1}   & &  & a^d_1 \\
 \rule[0mm]{0.1mm}{10mm} &  &  &  \rule[0mm]{0.1mm}{10mm} &   \rule[0mm]{0.1mm}{10mm} &  & & \rule[0mm]{0.1mm}{10mm} & &       \rule[0mm]{0.1mm}{10mm} & & &       \rule[0mm]{0.1mm}{10mm}  \\
 0 &0 & \rule[0.5ex]{6mm}{0.1mm} & 0 & a_{b}^1 &  a^1_{b-1}& \rule[0.5ex]{6mm}{0.1mm} & a_{b-d}^1 &  \cdots & a^d_{b} & a_{b-1}^d &  \rule[0.5ex]{6mm}{0.1mm} &  a^d_{b-d}\\
 \end{pmatrix}
$$
\vskip2mm
\noindent
By Corollary \ref{cor : Bezout bound in case B=R=1 }  any non-zero minor determinant of size $m\times m $   of $M$ will provide
a bound for $\cal Z_d(f)$ on $D_{\frac{1}{4}}$, since in Corollary \ref{cor : Bezout bound in case B=R=1 }, $\delta$ is the maximum of all non-zero minor determinants of size $m\times m $ of $M$.
But such a minor has to contain the first $d+1$ lines of $M$, as well as the last line. One sees that the absolute value $\Delta$ of any
non-zero determinant  of $ (d^2+d)\times(d^2+d)$ minor of the following
matrix
\begin{equation}\label{eq: the matrix tilde M}
\tilde{M}= \begin{pmatrix}
a^1_{d+1}   & \rule[0.5ex]{6mm}{0.1mm}   & a^1_1 & \cdots & a^d_{d+1}   & & \rule[0.5ex]{6mm}{0.1mm}  & a^d_1                           \\
 \rule[0mm]{0.1mm}{10mm}  &  &    \rule[0mm]{0.1mm}{10mm} &  & \rule[0mm]{0.1mm}{10mm} &  &  &    \rule[0mm]{0.1mm}{10mm}  \\
a^1_b  & \rule[0.5ex]{6mm}{0.1mm} & a^1_{b-d}& \cdots & a^d_{b} &  & \rule[0.5ex]{6mm}{0.1mm} &  a^d_{b-d}\\
 \end{pmatrix}
\end{equation}
\vskip2mm
\noindent
provides on $D_{\frac{1}{4}}$ the B\'ezout bound
\begin{equation}\label{eq: general Bezout bound }
 \cal Z_d(f)\le 5b \log(4+2(b+1) \frac{e^{2(d+1)^2\log(d+1)}}{\Delta}).
\end{equation}
$$ $$

\subsection{B\'ezout bound through the transcendence index}\label{Sec: Bezout bound through the transcendency index}
We start this subsection by the following definition of a notion of measure of the local transcendence of an analytic  transcendental function.

\bd\label{def:Transc.Ind}
For a transcendental analytic function $f: D_1\to \CC$ and for any $d\ge 1$ the
{\sl $d$-th transcendence index $\nu_d$ of $f$} is the maximal (with respect to all non-zero polynomials $P_d\in  \cal P_d$) multiplicity at $0$ of the function $g(z)=P_d(z,f(z))$. The non-decreasing sequence $\nu(f)=(\nu_1,\ldots,\nu_d,\ldots)$ is called the
{\sl transcendence sequence  of $f$}.
\ed

\br
The d-th transcendence index   measures the maximal order of contact at the origin between the graph of $f$ and algebraic curves of degree at most $d$. The higher this index is,
the less $f$ seems transcendental, since infinite $\nu_d$ means
that $f$ is algebraic.

\er

\br\label{rmk: bound from below of the transcendency index nu} As already observed in Proposition \ref{prop:max multiplicity and Bautin index}, $\nu_d$ is the Bautin index of the linear family associated to $f$ and the monomials $X^iY^j$ of degree $i+j$ at most $d$. Since the number of monomials in two variables of degree at most $d$ is $(d+1)(d+2)/2$, by Remark \ref{rem:transcendeal functions and Bautin index},
$\nu_d\ge (d^2+3d)/2$ always.
\er

\br\label{rmk: double bound between b and nu}
Using the notation $b$, as in the beginning of section \ref{Sec: Bezout bounds for transcendental analytic curves}, for the Bautin index of the family
$ \sum_{0\le j\le d} p_j(z)f^j(z)$, with $\deg p_j\le d$,
one has by Proposition \ref{prop:max multiplicity and Bautin index}
that $b= \mu$, where $\mu $ is the maximal multiplicity, with respect to the coefficients of the polynomials  $p_j$, of this family. Therefore
 one gets $\nu_d\le b= \mu \le \nu_{2d}$.
 In particular, inequality
 (\ref{eq: general Bezout bound }) gives the following proposition.
\er
 \bp\label{pr: Bezout and index of transcendency}
Let $f:D_1\to \CC$ be a transcendental analytic function with
transcendence sequence  $(\nu_d)_{d\ge 1}$, then
on $D_{\frac{1}{4}}$
$$  \cal Z_d(f)\le 5 \nu_{2d} \log(4+2( \nu_{2d}+1) \frac{e^{2(d+1)^2\log(d+1)}}{\Delta}), $$
  where $\Delta$ is the absolute value of a certain non-zero minor determinant
of size $(d^2+d)\times(d^2+d)$ of the matrix $\tilde{M}$ defined
  at (\ref{eq: the matrix tilde M}).
\ep
Now we assume that $  f(z)=\sum_{k=0}^\infty a_kz^k$ has rational Taylor coefficients
$ a_k=\frac{m_k}{p_k}$, with the greatest common divisor $m_k\wedge p_k$ of $m_k,p_k$ equal to $1$ and $p_k>0$,  and let us introduce the following notation.
\bn\label{not: h}
For any $l\ge 1$, let $h_l=\max\{p_k; k=1, \ldots, l\}.$
\en
Under this assumption we can bound from below the non-zero determinant of the $(d^2+d)\times (d^2+d)$ minors of $\tilde{M}$,  in terms of the transcendence index $\nu_{2d}$ and the height bound $h_{\nu_{2d}}$.

\bp\label{prop:rat.coef.delta}
For $f$ as above and for $\Delta$ the absolute value of a non-zero minor determinant of size $(d^2+d)\times (d^2+d)$   of $\tilde{M}$ we have
$$
\Delta \ge  h_{\nu_{2d}}^{-d^2(d+1)\nu_{2d}}.
$$
\ep

\bpr
For simplicity write $h$ for $h_{\nu_{2d}}$ and $\nu$ for $\nu_{2d}$.
Using the notation of \ref{no: coefficients a_j^i},
write the coefficients $a^1_k$, $k\le \nu$, as rational numbers having the same denominator $D$, this common denominator being at most $h^\nu$.
Next we write the coefficients $a^j_i$, for $1\le j\le d, \ 1\le i\le \nu$, as rational numbers having the same denominator $D^d.$ Indeed, $a^j_i$ are sums of the products $a^1_{i_1}\cdots a^1_{i_j}, i_1+\cdots +i_j=i$, with  $j\le d$.
 Therefore the determinant of a  $(d^2+d) \times (d^2+d)$ minor of $\tilde{M}$
can be written as a rational number having for denominator
$$
D^{d^2(d+1)} \leq h^{d^2(d+1)\nu}.
$$
But then such a non-zero determinant cannot be smaller in absolute value
than $h^{-d^2(d+1)\nu}$.
\epr

An important special case is when  there exist polynomials $R(d)$ and $S(d)$, with positive coefficients, such that
\be
\label{eq:growth.cond}
\nu_d\leq R(d), \ h_l\leq e^{S(l)}, \ d, l\ge 1
\ee
Under this condition we can guarantee that   $\cal Z_d(f)$  grows at most polynomially in $d$.
\bt\label{thm: Bound through the transcendency index sequence}
Assume that $f$ has rational Taylor coefficients at the origin,
and that the growth conditions (\ref{eq:growth.cond}) are satisfied.
Then on $D_{\frac{1}{4}}$
$$\cal Z_d(f) \le T(d),$$
for a certain polynomial $T$.
\et
\bpr
Under condition (\ref{eq:growth.cond}) and
by Proposition \ref{prop:rat.coef.delta}, since $S$ is an increasing function, we have
$$
\frac{1}{ \Delta } \le
 h_{\nu_{2d}}^{d^2(d+1)\nu_{2d}}
 \le e^{S(R(2d)) d^2(d+1)R(2d)}=e^{U(d)}.
$$
Now by Proposition  \ref{pr: Bezout and index of transcendency},
on $D_{\frac{1}{4}}$, we easily have for instance
$$
\cal Z_d(f) \le
5 R(2d) \log (4+ 2(R(2d)+1) e^{2(d+1)^2\log(d+1)}e^{U(d)})$$
$$
\le 5 R(2d) \log (4R(2d)e^{2(d+1)^3}e^{U(d)})
$$
$$
\le 10 R^2(2d)+ 10R(2d)(2(d+1)^3+U(d))
$$
\epr

\br\label{remark:transc.seq}
Producing instances of Taylor series $ f(z)=\sum_{k\ge 0} a_kz^k$ converging on $D_1$, with rational coefficients $a_k$ having denomi\-na\-tors boun\-ded from above by $e^{S(k)}$, where $S$ is a certain polynomial, is easy. Nevertheless the second assumption of Theorem \ref{thm: Bound through the transcendency index sequence}, concerning the growth of the transcendence sequence of $f$, is more difficult to control. A polynomial bound $\nu_d(f)\leq R(d)$ is known for solutions of some classes of algebraic ODE's (see \cite{Bin2016, Ga1999, Ne2009, Ne1986}). We expect such a bound to hold for Taylor series produced by some natural classes of recurrence relations. In Section \ref{Sec:lacun.} we give conditions on lacunarity of the series $f$ that allow estimates of the growth of $\nu_d(f)$. In general, we consider bounding of the growth of the transcendence sequence of $f$ as an important open question.

If we consider polynomials $P(z,y)=p_1(z)y+p_0(z)$ of degree $1$ in $y$, with $p_0(z),p_1(z)$ of degree $d$ in $z$, we are in the framework of the classical Pad\'e approximation. In this case the sequence of maximal multiplicities $\mu_d=\mu_d(f)$ of $g(z)=P(z,f(z))=p_1(z)f(z)+p_0(z)$ has the following remarkable description (see, for instance, \cite{Nik.Sor}): let
$$
f(z)=\frac{1}{q_1(z)+\frac{1}{q_2(z)+...}}, \ \ \deg q_l=s_l, \ l=1,2,...,
$$
be a continued fraction representation of the series $f(z)$. Then $\mu_d= s_1+s_2+\cdots+s_d$.
\er
For polynomials $P(z,y)=p_0(z)$ of degree $0$ in $y$ the behavior of $\mu_d$ was studied in \cite{Fri.Yom}, in particular, it was related there to linear non-autonomous recurrence relations for the Taylor coefficients of $f$.

\subsection{B\'ezout bound through the  Bautin multiplicity}\label{Sec: Bezout  Bound through the Bautin multiplicity}
%
%
In the previous section \ref{Sec: Bezout bound through the transcendency index}, in Proposition \ref{prop:rat.coef.delta}, some minors of the matrix $\tilde{M}$ were bounded from below, in terms of the transcendence index $\nu_{2d}$   and the height bound $h_{\nu_{2d}}$. On this base, on $D_{\frac{1}{4}}$,   $\cal Z_d(f)$  was  bounded from above by a polynomial in $d$ (see Theorem \ref{thm: Bound through the transcendency index sequence}),
under the condition that the sequences $(\nu_d)_{d\ge 1}$ and $(\log h_d)_{d\ge 1}$ are polynomially bounded.

A special case in which the transcendence index   (or thanks to the double inequality of Remark \ref{rmk: double bound between b and nu}, in which the Bautin index $b$ itself of the family
$ \sum_{j=0}^dp_jf^j, \deg p_j\le d$) is
polynomially bounded is the case that $b$ is minimal, that is equal to $d^2+2d$.
In this case the matrix $\tilde{M}$ of (\ref{eq: the matrix tilde M}) is an invertible square matrix of size $ d^2+d$ with the same determinant as the Bautin matrix $M$.

\bn
Being zero or not, let us call this $(d^2+d)\times(d^2+d)$ determinant {\sl the Bautin determinant}  of the family $ \sum_{j=0}^dp_jf^j, \deg p_j\le d$, and let us denote it by $\Delta_d$.
\en

So far the study has been done by looking at Taylor series at the origin. Let us now allow  Taylor expansions of $f$ at points $z$ near the origin. In this situation, the Bautin matrix $M$, as well as its submatrix $\tilde{M}$, have entries $a_i^j(z)$
that are analytic functions in the variable $z$. To emphasize this dependency, we adopt the notation  $\tilde{M}(z)$ and
$\Delta_d(z)$ (and keep the notation $\tilde{M}$ and $\Delta_d$ for
$z=0$).
 We shall assume that for each degree $d$ the Bautin determinant $\Delta_d(z)$, as a function of $z$, does not vanish identically. This is in particular true when $f(z)$ is a hypertranscendental function, i.e. all the derivatives are algebraically independent.
In this situation, for a generic base point $z$, $\Delta_d(z)\not=0$,  $\nu_d(z)\le d^2+2d$,
and therefore the transcendence index $\nu_d(z)$ is polynomially
bounded. The study of section \ref{Sec: Bezout bound through the transcendency index} could then be done by shifting the origin at some generic point $z$, however in this
translation one loses the control on the rationality of the coefficients of the Taylor
expansion of $f$, an assumption that is necessary to formulate at some fixed point (namely the origin, for simplicity), since this assumption makes no sense at generic points $z$. Nevertheless, still in the case that $\Delta_d(z)$ does not vanish identically, one can use the following dichotomy:
\vskip1mm
- when $\Delta_d(0)\not=0$, as just observed, we are in particular in the frame of Section \ref{Sec: Bezout bound through the transcendency index} where $\nu_d$ is polynomially bounded (by $d^2+2d$),
\vskip1mm
- when $\Delta_d(0)=0$, we can expand $\Delta_d(z)$ as a non-zero Taylor series in $z$ at the origin and study the multiplicity of this expansion with respect to $d$. We will see in Theorem \ref{thm: Bezout  Bound through polynomial hypertranscendence} that when the sequence of these multiplicities
has at most a polynomial growth, then a B\'ezout bound for $f$ is still possible.
\vskip1mm
\br
This dichotomy means that when some transversality defect for $f$ is quantitatively well controlled (through the multiplicity of $\Delta_d(z)$), then a good zero-counting bound
is possible, and for instance, in turn, a good bound for the density
of rational points of bounded height in the graph of $f$ will also be possible (see Theorem \ref{thm: bound for rational points of bounded height }).
\er
\bd\label{def:Hyper.Transc.Ind}
For any $d\ge 1$,  the {\sl $d$-th Bautin multiplicity $\eta_d$ of $f$} is the multiplicity at $0$ of the Bautin determinant $\Delta_d(z)$, considered as an analytic function of $z$.
The sequence $\eta(f)=(\eta_d)_{d\ge 1}$ is called the {\sl Bautin multiplicity sequence of $f$}.
\ed
In brief, for each $d\ge 1$, we can write
$$
\Delta_d(z) = \alpha_d z^{\eta_d}+O(z^{\eta_d+1}), \hbox{ with } \alpha_d\ne 0.
$$


\bt\label{thm: Bezout Bound through polynomial hypertranscendence}
Assume that $f:D_1\to \CC$ is an analytic function with
rational Taylor coefficients at the origin satisfying  the growth condition (\ref{eq:growth.cond}) and such that
there exists a polynomial $R$ with
$\eta_d\le R(d)$, for $d\ge 1$.
Then on $D_{\frac{1}{4}}$
$$\cal Z_d(f) \le T(d),$$
for a certain polynomial $T$.
\et

\bpr
\label{rem: eta bounds nu}
The Bautin matrix
of Definition \ref{de: Bautin matrix} is the matrix of size $(d+1)^2$, in the base of monomials of
${\cal Q}_d$,
of the linear map $\widehat{\Lambda}: {\cal Q}_d \to (\CC\{z\})^{(d+1)^2}$ sending a polynomial
to the vector of the first $ (d+1)^2$ derivatives
$(P(z,f(z))^{(j)}/j!$ of $P(z,f(z))$. For a given polynomial $P\in {\cal Q}_d$,
such that the multiplicity at the origin of $P(z,f(z))$ is maximal, and therefore is the Bautin index $b_d$
of $f$ for the family of monomials of ${\cal Q}_d$, one can write $b_d=d^2+2d+r$, with $r\ge 0$. Then the multiplicity
at the origin of the first $(d+1)^2$ derivatives of $P(z,f(z))$ is bigger than $r$. Now writing the Bautin matrix in a basis of ${\cal Q}_d$ starting with $P$, the elements of the first column of this matrix have multiplicity at least $r$. It follows that the Bautin determinant $\Delta_d(z)$ itself has multiplicity at least
$r=b_d-d^2-2d$, and therefore $\eta_d\ge b_d-d^2-2d\ge \nu_d-d^2-2d$.
As a conclusion, when $\eta_d$ is polynomially bounded, $\nu_d$ is polynomially bounded as well.
The existence of the polynomial $T$ then follows
from Theorem \ref{thm: Bound through the transcendency index sequence}.
\epr

\br
\label{rem: arity and degree of Delta}
The function  $\Delta_d(z)$ is a polynomial with
coefficients in $\ZZ$ in the variables
$a_1(z), \ldots, a_{d^2+d}(z)$ and with degree
$\frac{d(d+1)^2}{2}$. Indeed, the entries of $\tilde{M}$
are the functions $a_i^j$, with $i=1,\ldots, d^2+d$, $j=1, \ldots, d$, and $a_i^j$ is a sum of products of type $a_{i_1}\cdots a_{i_j}$, $i_1+\cdots +i_j=i$. The Bautin multiplicity $\eta_d$ therefore measures the degree
of cancellation allowed by the polynomial $\Delta_d$ applied on the $d^2+d$ first derivatives of $f$.
This simple observation suggests to introduce,
for an hypertranscendental analytic function $f: D_1\to \CC$, the notion of
{\sl polynomial  hypertranscendence}, defined by  the existence of polynomials $A,B\in \RR[X]$ with positive coefficients, such that
for any $d\in \NN$, for any polynomial $P\in \ZZ[X_0, \ldots, X_{A(d)}]$, with degree $\le  d$,
the multiplicity of $P(\frac{f(z)}{0!}, \frac{f'(z)}{1!}, \ldots, \frac{f^{(A(d))}(z)}{A(d)!})$ at the origin is bounded by $B(d)$.
This notion of strong hypertranscendence, relevant by itself, is motivated here by the remark that since there exists $p\in \NN$ such that $d^2+d \le d^p$ and $\frac{d(d+1)^2}{2}  \le A(d^p)$, we have $\eta_d\le B(d^p)$ for such a function.
\er


\subsection{Lacunary series}\label{Sec:lacun.}

The aim of this section is to give concrete instances of series $f$ satisfying polynomial growth condition for the transcendence sequence $(\nu_d)_{d\ge 1}$  considered in (\ref{eq:growth.cond}), and thus having a polynomial B\'ezout bound. For this we focus on lacunary series for which the computation of a bound for the transcendence sequence $(\nu_d)_{d\ge 1}$  is possible.
This case has been considered in \cite[Theorem 6.1]{Com.Pol}, where a family of analytic functions with $\cal Z_d(f)$ having prescribed growth is given. Our conditions, being more flexible, improve on these earlier
conditions (see Remark \ref{rm: Comparaison Lacunarity}).


We begin with the following remark which improves on the estimates of Propositions \ref{prop:rat.coef.delta} and \ref{thm: Bound through the transcendency index sequence}, in case of lacunary series.
%
\br\label{rm: density of nonzero Taylor coefficients} Assume that the Taylor coefficients
of the series $f$ are rational numbers and denote by $\theta_d$ the amount of those non-zero coefficients among
$a_0, \ldots, a_{\nu_d}$. Then, with
exactly the same proofs adapted to this notation, Proposition \ref{prop:rat.coef.delta} and Theorem \ref{thm: Bound through the transcendency index sequence} may be formulated as follows:
the absolute value $\Delta$ of a non-zero minor determinant
of size $(d^2+d)\times (d^2+d)$   of $\tilde{M}$ satisfies
$ \Delta\ge h_{\theta_{2d}}^{-d^2(d+1)\theta_{2d}}$,
and consequently when there exist polynomials $R $ and $S$ such that
$ \theta_d\le R(d), \ \ h_d\le e^{S(d)}$, one gets
$ \cal  Z_d(f) \le T(d) \nu_{2d}\log \nu_{2d}, $
for some polynomial $T$.
\er


Now in order
to  estimate the growth of the transcendence sequence $(\nu_d)_{d\ge 1}$, let us assume that the lacunarity of the series $f$ is quantitatively controlled by the following condition
%
\begin{equation}\label{eq: squarre lacunarity condition }
 f(z)=\sum_{k=1}^\infty a_kz^{n_k},  \
a_k\ne 0 \ \hbox{ and for any }\  k\ge 1,  \ n_{k+1}>n_k^2.
\end{equation}
Note that in what follows no assumption is made
on the rationality of the Taylor coefficients of $f$. We assume that
$f$ is analytic on $D_1$ and bounded there by $1$.
%
\bl\label{le:lacun.powers}
Under condition (\ref{eq: squarre lacunarity condition }), for any
$l\ge 1$
and for any $m,j\in [0, n_{l+1}-1]$, the series $z^mf^j(z)$ contains the non-zero monomial $(a_{l+1})^jz^{jn_{l+1}+m},$ and no other monomials with degree in $]jn_{l+1}+m, n_{l+2}[$.
\el
\bpr
For $j=0$ the series $z^mf^j(z)$ is the monomial $z^m$ and thus
is of the required from. Now for $j\geq 1$, write $f(z)$ as the sum $f(z)=f_{l+1}(z)+\bar f_{l+1}(z),$ with
$ f_{l+1}(z)=\sum_{k=1}^{l+1} a_kz^{n_k},$ and $ \bar f_{l+1}(z)=\sum_{k=l+2}^\infty a_kz^{n_k}.$
The monomial $(a_{l+1})^jz^{jn_{l+1}}$ has the highest degree among the monomials of the series $f^j(z)$ which come from the terms in $f^j_{l+1}$.
On the other hand, all other monomials in $f^j(z)$ which come from the products of terms in $f_{l+1}(z)$ and in $\bar f_{l+1}(z),$ have degree at least $n_{l+2}$. The terms of the series $x^mf^j(z)$ are the terms of $f^j(z)$ shifted by $m$ and therefore for $m,j < n_{l+1}$ we have
$$
jn_{l+1}+m\leq (n_{l+1}-1)n_{l+1}+ n_{l+1}-1=n_{l+1}^2-1<n_{l+2},
$$
since by condition (\ref{eq: squarre lacunarity condition }), $n_{l+2}>n_{l+1}^2$.
\epr
%
\bp\label{prop: transcendency index and lacunarity }
Under condition (\ref{eq: squarre lacunarity condition }),
for any  $l\ge 1$ and for any $d$ in the interval $[n_l,n_{l+1}-1]$ we have
$$
n_{l+1}\leq \nu_d \leq n_{l+1}^2-1 < n_{l+2}.
$$
\ep
\bpr
First of all, for $l\ge 1$ and $d\ge n_l$, we have
$\nu_d\ge \nu_{n_l}\geq n_{l+1},$
since for the polynomial $P(z,y)=y-\sum_{k=1}^l a_kz^{n_k}$ of degree $n_l$, the function $ P(z,f(z))=\sum_{k=l+1}^\infty a_kz^{n_k}$ has multiplicity $n_{l+1}$ at the origin.

Let now $P_d(z,y)=p_d(z)y^d+\cdots+p_1(z)y+p_0(z)$ be a polynomial of degree $d\leq n_{l+1}-1$ and let us prove that the multiplicity of $P_d(z,f(z))$ at the origin is at most $n_{l+1}^2-1$; this will prove that $\nu_d\le n_{l+1}^2-1< n_{l+2}$.

Denote by $s\leq d$ the highest degree of $y$ in $P_d(z,y)$ for which the polynomial $p_s(z)$ is not identically zero, and let us write $p_s(z)=ax^r+bx^{r-1}+\cdots, \ r\leq d,$ with $a\ne 0$.
By Lemma \ref{le:lacun.powers}, the summand $ax^rf^s(z)$ in $P_d(z,f(z))$ contains the  monomial $v=a(a_{l+1})^sz^{sn_{l+1}+r}.$
Let us show that this monomial cannot cancel with any other monomial in $P_d(z,f(z))$, because, if it is the case, since
for  $s,r\leq d<n_{l+1}$, we
have
$$  sn_{l+1}+r\le (n_{l+1}-1)n_{l+1}+n_{l+1}-1 \le n_{l+1}^2
<n_{l+2}, $$
it will finish the proof.
As just noticed, since $sn_{l+1}+r<n_{l+2}$, the monomial $v$ can cancel only with the monomials coming from the truncated series $f_{l+1}(z)= \sum_{k=1}^{l+1} a_kz^{n_k}$ introduced in the proof of Lemma
\ref{le:lacun.powers}, since
$f(z)-f_{l+1}(z)= \sum_{k=l+2}^\infty a_kz^{n_k}$.
But on one hand $v$ cannot cancel with any monomial in $p_s(z)f_{l+1}^s(z)$, and on the other hand, for any $q<s$, the monomials in $f^q(z), \ q<s,$ coming from the truncated series $f_{l+1}(z)$, have   degree at most $qn_{l+1}$. Hence the highest degree of monomials in $p_q(z)f_{l+1}^q(z)$ can be $qn_{l+1}+d<(q+1)n_{l+1}\leq sn_{l+1}+r$. As announced, we
conclude that $v$ cannot cancel with any other monomial in $P_d(z,f(z))$.
\epr


Lemma \ref{le:lacun.powers} does not only let us
 bound the terms of the sequence $\nu$ like in Proposition
 \ref{prop: transcendency index and lacunarity },
it also lets us compute
some $(d+1)^2\times(d+1)^2$ non-zero minor determinant
in the Bautin matrix $M$ (up to allowing more than $b+1$ rows in $M$).

\bp\label{prop:Bautin.det}
Under condition (\ref{eq: squarre lacunarity condition }),
for any $l\ge 1$, for any
 $d\in [n_l,n_{l+1}-1]$, there exists a
 $(d+1)^2\times(d+1)^2$ minor in $M$ (up to allowing more than
 $b+1$ rows in the definition of the Bautin matrix $M$) with non-zero determinant $\Delta=(a_{l+1})^{\frac{1}{2}d(d+1)^2}.$
For $d=n_{l+1}-1$ this determinant
is the upper square $(d+1)^2\times (d+1)^2$ minor of the Bautin matrix $M$.
\ep


\bpr Let us fix $l\ge 1$, and $d\in [n_l, n_{l+1}-1]$. Then
for any $j=0,\ldots,d$ by Lemma \ref{le:lacun.powers},
we have in $M$
a lower-triangular square $(d+1)\times (d+1)$ block $M_j$
corresponding to
lines ranging from  $jn_{l+1}$ to $jn_{l+1}+d$ and columns
ranging from $j(d+1)$ to $(j+1)(d+1)-1$ in $M$,
and $M_j$ has for entries $(a_{l+1})^j$ on its main diagonal. Note that
$dn_{l+1}+d$ may be bigger than $b+1$, so we maybe have to consider
a matrix having more lines then the Bautin matrix $M$, which is harmless.

Now, if we drop from $M$ all the lines which are not in $M_j$ for some $j=0,\ldots, d,$ we obtain a $(d+1)^2\times (d+1)^2$ minor
$M'$ which is lower triangular and has the blocks $M_j, \ j=0,\ldots, d,$ on its main diagonal.
Hence its determinant $\Delta$ is
$(a_{l+1})^{\frac{1}{2}d(d+1)^2}.$

For $d$ equal to its maximal value $n_{l+1}-1$ each line of $M$ belongs to one block $M_j$, and hence $M'$ coincides with the upper square $(d+1)^2\times (d+1)^2$ minor of the Bautin matrix $M$.
\epr

So far, in condition (\ref{eq: squarre lacunarity condition })
we have required that the lacunarity of the series $f(z)=\sum_{k=1}^\infty a_kz^{n_k} $ is big
enough.
We now require in addition that the lacunarity of $f$ is not too big, in the following condition
\begin{equation}\label{eq: lacunarity bounded from above}
\hbox{ There exists } q>2, \ \hbox{ such that for any }
k\ge 1, \ n_k^2< n_{k+1} \le n_k^q.
\end{equation}
Under this assumption we can now
show that $f$ has a polynomial B\'ezout bound on
$ D_{\frac{1}{4}}$.

\bt\label{thm: Bezout Bound for lacunary series}
Under condition (\ref{eq: lacunarity bounded from above}),  the transcendence sequence $\nu_d$ grows at most polynomially in $d$. More accurately we have
$$
\nu_d(f) <  d^{q^2}.
$$
If, in addition, for a certain fixed $p>0,$ and for any
$k\ge 1$
we have $|a_k|\geq e^{-n_k^p}$, then
on $ D_{\frac{1}{4}}$
$$
\cal Z_d(f) \le  10(2d)^{q^2} ( 1 +qd^2+ 5d^{pq+3} ).
$$
\et
\bpr
Let the degree $d\ge 1$ be given and let $l$ be such that $d\in [n_l,n_{l+1}-1]$. By Proposition \ref{prop: transcendency index and lacunarity } we have
$$
\nu_d < n_{l+2} \leq (n_l)^{q^2} \leq d^{q^2}.
$$
Now by Proposition \ref{prop:Bautin.det} we obtain the existence
of some $(d+1)^2\times(d+1)^2$ minor in $M$ with non-zero determinant $\Delta$, such that
$$
\vert \Delta \vert = \vert a_{l+1}\vert^{\frac{1}{2}d(d+1)^2}
\geq \exp (-\frac{1}{2}d(d+1)^2n_{l+1}^p)$$
$$\geq  \exp(-\frac{1}{2}d(d+1)^2n_{l}^{qp})
\geq \exp(-\frac{1}{2}d^{qp+1}(d+1)^2).
$$
Now by Proposition \ref{pr: Bezout and index of transcendency},
$$
\cal Z_d(f) \le  5
\nu_{2d}
\log ( 4+\frac{\nu_{2d}+1}{\Delta}e^{2(d+1)^3} )\leq
$$
$$
\leq (2d)^{q^2}5 \log ( 4+((2d)^{q^2}+1) e^{\frac{1}{2}d^{qp+1}(d+1)^2+2(d+1)^3} ).$$
Since for instance,
$1+(2d)^{q^2}\le e^{2dq^2}$, $\frac{1}{2}d^{qp+1}(d+1)^2+2(d+1)^3 \le 5d^{pq+1}(d+1)^2
\le  10d^{pq+3} $, we have
$$\cal Z_d(f) \le  (2d)^{q^2}5
\log ( 4+ e^{2qd^2+10d^{pq+3}} ).  $$
And finally since $\log(4+e^x)\le 2+x$, for $x\ge 0$, we obtain
$$ \cal Z_d(f) \le  10(2d)^{q^2} ( 1 +qd^2+ 5d^{pq+3} ). $$
This completes the proof of Theorem \ref{thm: Bezout Bound for lacunary series}
\epr

\br\label{rmk: special degrees}
 When the gaps between the degrees $n_l$ grow faster than assumed in Theorem \ref{thm: Bezout Bound for lacunary series}, that is faster than forced by condition (\ref{eq: lacunarity bounded from above}), the asymptotic growth of the transcendence indices $\nu_d$ and of bounds on $\cal Z_d(f)$ fastens accordingly. Notice, however, that from Proposition
 \ref{prop: transcendency index and lacunarity },
 for a subsequence of degrees of the form $d=n_{l+1}-1, \ l\ge 1$,
 under condition (\ref{eq: squarre lacunarity condition }) only,  the multiplicity $\nu_d$ is at most $d^2+2d$. Under appropriate assumptions on the coefficients $a_l$ for the above subsequence of degrees $d$ we obtain a polynomial bound also for $\cal Z_d(f)$.
 This phenomenon is to compare to a similar behaviour
 in  \cite[Theorem 1.1]{Com.Pol2007} or \cite[Theorem 0.3]{Su2002}, \cite[Theorem 1.3]{Su2006}, where the minimal asymptotic
 for $\cal Z_d(f)$ is obtained on some sequence of degrees $d$ going to infinity.
\er

\br\label{rm: Comparaison Lacunarity}
In \cite[Corollary 6.2]{Com.Pol}, for $\alpha\ge 3$, the lacunarity condition,
\begin{equation}\label{eq: Coman and Poletsky lacunarity condition}
n_{k+1}=n_k^{\alpha-1}, \ \ a_k=e^{-n_k\log n_k}
\end{equation}
is considered as a condition implying a polynomial B\'ezout bound for $f$ (with prescribed
order of asymptotic in $[\alpha-1, \alpha]$ as explained in the introduction).
Here our condition (\ref{eq: lacunarity bounded from above}) and our assumption
$\vert a_k\vert \ge e^{-n_k^p}$, for some $p>0$,
in Theorem \ref{thm: Bezout Bound for lacunary series}, are somewhat more flexible as conditions
giving a polynomial B\'ezout bound for $f$.

\er

\subsection{ Rational Taylor coefficients via recurrence relations}\label{Sec:Rec.rel.}


In this section we study recurrence relations for the Taylor coefficients of the series $f(z)= \sum_{k=1}^\infty a_k z^k$, assuming that the starting coefficients are rational numbers. In case the recurrence relation is linear (with polynomial coefficients in $1/k$) it turns out that the $a_k$'s satisfy the bound $h_l\le e^{S(l)}, \ l\ge 1$
of  condition  (\ref{eq:growth.cond}), one of the two hypotheses  required in  Theorems \ref{thm: Bound through the transcendency index sequence}   and \ref{thm: Bezout  Bound through polynomial hypertranscendence}. Furthermore in case
$f$ is a solution of an algebraic dif\-ferential equation  with polynomial coefficients, the other hypothesis required in Theorem \ref{thm: Bound through the transcendency index sequence}, namely a bound that is  polynomial in $d$ for
the transcendence index $\nu_d$, is automatically satisfied (see Remark \ref{rmk:ODE} and Theorem \ref{thm:ODE}).

Let
$$
Q(k,u_1,u_2,\ldots,u_r)=\sum_{|\beta|\leq d_1}p_\beta(k)u^\beta
$$
be a polynomial of degree $d_1$ in the variables $u_1,u_2,\ldots,u_r$ and  with  coefficients  $p_\beta(k)=\sum_{i=0}^{d_2}c_{\beta,i}(\frac{1}{k})^i$ being polynomials in $\frac{1}{k}$ of degree $d_2$. We consider a polynomial recurrence relation of length $r$ of the form
%
\be\label{eq:rec.rel}
a_{k+1}=Q(k,a_{k},a_{k-1},\ldots,a_{k-r+1}).
\ee
We assume that the coefficients $c_{\beta,i}$ are rational numbers, as well as the initial terms $a_0,a_1,\ldots,a_{r-1}$ of the sequence $a=(a_0,a_1,\ldots,a_{r-1},a_r,\ldots).$ For $k\ge r-1$,
let
$D_k$ denote the common denominator of $a_0,\ldots,a_k$, when
those rational numbers are   written in their irreducible form.
We also denote by  $L_1$ (respectively, $L_2$) the common denominator of all the coefficients $c_{\beta,i}$, $i=0, \ldots, d_2 $, $\vert \beta\vert \le d_1$ (respectively, the common denominator of all the initial given terms $a_0,a_1,\ldots,a_{r-1}$) again when those rational numbers are written in their irreducible form. Note that $D_{r-1}=L_2$.
%
\bp\label{prop:denom.growth.rec.}
With the notation above, for any  $k\ge r-1$,
$$
D_k\leq e^{Md_1^{k-r+1}k\log k},
$$
where
$$
M=\max \left (\frac{\log L_2}{(r-1)\log(r-1)}, \ d_2+\frac{\log L_1}{\log 2} \right ).
$$
\ep
\bpr
The products of $a_j$ entering $Q$ in (\ref{eq:rec.rel}) can be written with denominator $D_k^{d_1}$. Therefore, the next term $a_{k+1},$ and hence all the terms $a_0,\ldots,a_k,a_{k+1},$ can be written with the common denominator $\tilde D_{k+1}=L_1k^{d_2}D_k^{d_1}.$ Now we prove by induction that

\be\label{eq:induction.ineq}
D_k\leq e^{Md_1^{k-r+1}k\log k}.
\ee
For $k=r-1$ we have $D_{r-1}=L_2,$ and (\ref{eq:induction.ineq}) is satisfied by the choice of $M$. Assuming that the required inequality is satisfied for a certain $k\geq r-1$, we now prove it for $k+1$. We have
$$
D_{k+1} \leq \tilde D_{k+1}=
L_1k^{d_2}D_k^{d_1} \leq L_1k^{d_2} e^{Md_1^{k+1-r+1}k\log k}
$$
$$
= e^{Md_1^{k+1-r+1}k\log k+\log L_1 + d_2\log k}
$$
$$
=e^{Md_1^{k+1-r+1}
(k\log k+\frac{\log L_1 + d_2\log k}{Md_1^{k+1-r+1}})}.
$$
By the choice of $M$ the last expression does not exceed
$$
e^{Md_1^{k+1-r+1}(k\log k+\log k)}<
e^{Md_1^{k+1-r+1}(k+1)\log (k+1)}.
$$
This completes the proof of Proposition \ref{prop:denom.growth.rec.}.
\epr

Notice that for $d_1>1$ the denominators grow as a double exponent, i.e. faster than an exponent of a polynomial in $k$. The trivial  example of recurrence relation $a_{k+1}=a_k^2$, i.e. $a_k=a_0^{2^k},$ shows that this growth indeed happens in recurrence relations of the form (\ref{eq:rec.rel}).
%
\br\label{rmk:ODE}
However, in the special case of linear recurrence relations of the form (\ref{eq:rec.rel}), we have $d_1=1$, and the bound of  Proposition \ref{prop:denom.growth.rec.} takes the form $D_k\leq e^{Mk\log k}$. This special case includes Poincar\'e-type recurrence relations, which are satisfied by the Taylor coefficients of solutions $f(z)$ of linear differential equations with polynomial coefficients.
In the more general case where $f$ satisfies an algebraic differential equation
$f^{(d)}=Q(z,f(z), \ldots, f^{(d-1)})$, where $Q$ is some given polynomial in $\QQ[X_1,\ldots, X_d]$,
iteration of derivation of each member of this equation leads to equations of type
$f^{(k)}(z)=Q_k(z,f(z), \ldots, f^{(d-1)})$, where $Q_k$ is a polynomial in
$\QQ[X_1,\ldots, X_d]$ with  controlled height of its
coefficients and controlled degree with respect to $d$. Studying these derivations and using some results of \cite{BinarXiv} one obtains also in
this case the bound on the height of $a_k=f^{(k)}(0)$ required by our growth condition   \eqref{eq:growth.cond}.
Therefore, combining well-known bounds on the transcendence sequences $(\nu_d)_{d\ge 1}$ of solutions of differential equations with polynomial coefficients, that turn out to be polynomially bounded in $d$ (see
\cite{Bin2016, Ga1999, Ne1986, Ne2009}), our Proposition \ref{prop:denom.growth.rec.} and Theorem \ref{thm: Bound through the transcendency index sequence} above, we immediately
obtain the following statement,
 a result, which was recently proved (among others results in this direction and by others methods) in \cite[Corollary 4, Theorem 6]{BinarXiv} (note that in \cite{BinarXiv}, no assumption on the rationality of initial conditions is required).

\er
%
\bt\label{thm:ODE} Let $f$ be an analytic function,  defined on the unit disc, that is solution of an algebraic differential equation
 with rational coefficients and initial conditions.
 Then there exists a polynomial  $T$
 such that on $D_{\frac{1}{4}}$, $\cal Z_d(f)\le T(d)$.
\et

%
\section{ Bautin determinant for random series}\label{Sec: Bautin determinant for random series}
\setcounter{equation}{0}

In this section, we discuss the behaviour of the Bautin determinant for random Taylor coefficients. We prove
that for any $p\in ]0,1[$, there exists a set $E_p$ of probability $p$, such that for any series $f\in E_p$, the corresponding
Bautin determinant, as a function of $d$,
is bounded from below by
$e^{U_p(d)}$, for a certain polynomial $U_p$. In case $\Delta_d\not= 0$, the Bautin index of the family $\sum_{j=0}^dp_jf^j$ is $d^2+2d$, and consequently the transcendence index $\nu_d$ is bounded by $d^2+2d$, and thus the growth condition on the transcendence indices sequence in Theorem
\ref{thm: Bound through the transcendency index sequence}  is fulfilled. It follows that for any $p\in ]0,1[$, for Taylor coefficients in a set of probability $p$, $\cal Z_d(f)$ is polynomially bounded in $d$.

Let us fix some integer $d\ge 1$ and let us start by the following remark.

\br \label{rmk: arity and degree of the Bautin determinant}
As already noticed in Remark \ref{rem: arity and degree of Delta}, as a polynomial in the Taylor coefficients of the series $f$, the Bautin determinant (of size $d^2+d$), still denoted $\Delta_d$, is a polynomial in
the variables
$a_1, \ldots, a_{d^2+d}$
 and with degree
$ \frac{d(d+1)^2}{2}$.
\er
For $d\ge 1$, following this remark, and to be more general, we will consider instead of the polynomial $\Delta_d$ of arity
(number of variables)
 $d^2+d$ and degree $d(d+1)^2/2$, any polynomial
with arity and degree polynomially bounded
in $d$.

%



Let $I=[-1, 1] \subset \RR$.
 We consider the unit $n$-dimensional cubes
 $I^n\in \RR^n, \ n\ge 1$,  and the infinite dimensional unit cube $I^\infty = \varprojlim_{n\in \NN}I^n$ that comes with its standard projections $\pi_n : I^\infty \to I^n$. Let us denote by $\mu_n$ the probability Lebesgue measure on $I^n$, for any $n\ge 1$.
For any $n\ge 1$ and any measurable set $G\subset I^n$
 denote by $\tilde G\subset I^\infty$ the cylinder $\pi_n^{-1}(G)$ over $G$.
 The probability measure $\mu$ on $I^\infty$ is defined by setting $\mu(E)=\sum_{i=1}^{\infty} \mu_{n_i}(G_{n_i}),$
 for any subset
 $E \subset I^\infty$  that can be expressed as a disjoint union
 of cylinders $\tilde{G}_{n_i}$, with $G_{n_i}$ a $\mu_{n_i}$-measurable in $I^{n_i}$.



We identify the sequences $(a_k)_{k\ge 0}\in I^\infty$ with the analytic functions $f(z)=\sum_{k=0}^\infty a_kz^k$, this series converging at least in the interior of $D_1$. For a polynomial $Q$ of arity $m$ and for $f=(a_k)_{k\ge 0}\in I^\infty$ we denote $Q(a_0, \ldots, a_{m-1})=Q(\pi_{m}(f))$ by $ Q(f)$.



Let finally $(Q_d)_{d\ge 1}$ be a sequence of polynomials $Q_d$, of   arity  $m_d$ and   degree $q_d$, and let us assume
\begin{equation}\label{eq: Bautin determinant bigger than 1}
\vert Q_d \vert_{I^{m_d}}=\max \{ \vert Q_d(x)\vert ; x\in I^{m_d} \}\ge 1.
\end{equation}
 Note that we also have $ \vert \Delta_d \vert_{I^{d^2+d}}\ge 1$, since
 $\Delta_d=1$ when  for instance $a^1_{d+1}=1$ and $a_i^j=0$  for $i=1, \ldots, d^2+d$,
 $j=1, \ldots, d$, $i\not=d+1$, $j\not=1$.
%
\bt\label{thm: Bound for Delta at random f }
With the above notation, for any $p\in ]0,1[$, there exists a set $E_p\subset I^\infty$ of measure $p$, such that for any $f\in E_p$, and for any $d\ge 1$,
$$
|Q_d(f)|\geq \left (\frac{3(1-p)}{2 \pi^2 d^2 m_d}\right )^{q_d}.
$$
In particular, for $q_d,m_d$ satisfying $q_d\leq d^{\kappa_1},
m_d \leq d^{\kappa_2},$ for some $\kappa_1, \kappa_2\ge 0$,
we have, with probability at least $p$:
$$
|Q_d(f)|\ge  e^{-(\gamma_p  + \kappa_2)d^{\kappa_1+1} },
$$
where $\gamma_p$ goes to $+\infty$ as $p$ goes to $1$.
\et
%
\bpr Let $p\in ]0,1[$, $d\ge 1$,
and $\theta_d=\frac{6(1-p)}{\pi^2 d^2}$.
 We define the real number $\varepsilon_d$ as the maximum of the numbers  $\varepsilon$ such that the set
$$
V_d=\{u \in I^{m_d}, \ |Q_d(u)| \leq \varepsilon \}
$$
satisfies $\mu_{m_d}(V_d)\leq \theta_d$.
Now for $V= \bigcup_{d=0}^\infty V_d \subset I^\infty$ and $E_p=I^\infty\setminus V$ we have
$$\mu(V)\leq \sum_{d=0}^\infty \theta_d=1-p, \hbox{ and thus }\ \mu(E_p )\ge p. $$
Clearly, for any $f\in E_p$, for any $d\ge 1$,
we have $|Q_d(f)|\geq \varepsilon_d$.
It remains to estimate the numbers $\varepsilon_d$.
We use for this purpose the following multivariate Remez inequality proved in \cite{Bru.Gan} (see also \cite{Bru1999}).

\smallskip
\noindent {\it Let   $Z$
be a measurable subset of $I^n$.
Then for every real polynomial $P$ in $n$ variables and of degree $d$,
%
\be\label{eq:BG}
\vert P \vert_{I^n}  <
\left(\frac{4n}{\lambda}\right)^d
\vert P \vert_Z,
\ee
where $\lambda= \mu_n(Z)$.}
%
%
\smallskip

\noindent
Applying  inequality (\ref{eq:BG}) to $Z=V_d$ and $P=Q_d$,
by our assumption (\ref{eq: Bautin determinant bigger than 1}) we get
$$
1 \leq |Q_d|_{I^{q_d}} < \left (\frac{4m_d}{\mu_{m_d}(V_d)}\right )^{q_d}\varepsilon_d,
$$
or equivalently
$\mu_{m_d}(V_d) \leq 4m_d\varepsilon_d^{\frac{1}{q_d}}.$
In particular, for $\varepsilon_d= \left (\frac{\theta_d}{4m_d}\right )^{q_d}= \left(\frac{3(1-p)}{2 \pi^2 d^2 m_d}\right )^{q_d}$, we have
$\mu_{m_d}(V_d)\leq \theta_d.$
This completes the proof of the first inequality of Theorem \ref{thm: Bound for Delta at random f }. Substituting into this inequality $q_d=d^{\kappa_1}, m_d=d^{\kappa_2},$ we obtain for any $f\in E_p$ and for any $d\geq 1$
$$ \vert Q_d(f) \vert \ge \left( \frac{3(1-p)}{2\pi^2} \right)^{d^{\kappa_1}}d^{-2-\kappa_2d^{\kappa_1}}
=e^{-c_p d^{\kappa_1}}e^{ -(2+\kappa_2d^{\kappa_1})\log d }$$
$$
\ge e^{-c_p d^{\kappa_1}}e^{ -(2+\kappa_2d^{\kappa_1})d } \ge e^{ -(2+c_p+\kappa_2)d^{\kappa_1+1} }=e^{ -(\gamma_p+\kappa_2)d^{\kappa_1+1} },
$$
where $ c_p=\log \left(\frac{2\pi^2}{3(1-p)}\right)>0, \ \gamma_p=2+c_p$.

\epr

We apply Theorem \ref{thm: Bound for Delta at random f } to the case where
$Q_d=\Delta_d$, then
as a consequence of Theorem \ref{thm: Bound through the transcendency index sequence}, we obtain the following statement.

%
\bc\label{cor: Bezout at random}
With the above notation, for any $p\in ]0,1[$, there exists a set $E_p\subset I^\infty$ of measure $p$, such that for any
$f\in E_p$, on $D_{\frac{1}{4}}$,
$\cal Z_d(f)\le C_pd^8$, where $C_p\to +\infty$ as $p\to 1$. Or, in other words, with probability $1$, random series satisfy  polynomial B\'ezout bounds (with degree at most $8$).
\ec

\section{Analytic functions with few rational points in their graph}\label{Sec:Rat.Points}
\setcounter{equation}{0}
We start this section by the following definition.
%
\bd Let  $x=(x_1, \ldots, x_n)\in \QQ^n$. The height of $x$ is the integer
$\max\{ \vert a_i \vert, \vert b_i \vert; i=1, \ldots, n \}$, where
$x_i=a_i/b_i$ with $a_i, b_i\in \ZZ$, $a_i\wedge b_i=1$, $i=1, \ldots, n$.
\ed
Explicit bounds on the number $\#X(\QQ,T)$ of rational points $X(\QQ,T)$ of height at most $T$,
in some given set $X\subset \RR^n$, are usually related to B\'ezout bounds satisfied
by $X$.
Let us  assume for instance that $X$ is a transcendental set definable in some o-minimal structure expanding the real field, and to simplify, of dimension $1$. Then following   \cite{PiWi2006}, that generalizes the classical by now method of \cite{BoPi1989}, one knows that
$X(\QQ,T)$ is contained in a certain number $H_{X,T,d}$ of hypersurfaces of $\RR^n$ of degree $d$, this number being of the form
$C_{X,d}T^{\tau_d}$, with $\tau_d\to 0$ when $d\to \infty$.
It follows that since the definable set $X$ satisfies a B\'ezout bound
(see \cite[Theorem 1]{BoPi1989}, \cite{PiWi2006})
%
\begin{equation}\label{eq: Pila Wilkie bound}
\forall \epsilon>0, \exists C_{X,\epsilon},  \forall T\ge 1,  \
\#X(\QQ,T)\le C_{X,\epsilon} T^\epsilon.
\end{equation}
Now in case the curve $X$ is given by a system of convenient parametrizations
(as  mild parametrizations defined in \cite{Pi2006}, or slow parametri\-za\-tions defined
in \cite{CoMi2016}), or more simply, in case $X$ is the graph $\Gamma_f$ of some transcendental analytic function $f$ on a compact interval of $\RR$, a computation shows
 that the constant $C_{\Gamma_f,d}$ is polynomially bounded in $d$ and
that $T^{\tau_{\lfloor \log T\rfloor}}$ is a constant $K$
independent of $T$ (see \cite[Proposition 2.4]{Pi2006}, \cite[Proposition 2.18]{CoMi2016}).
Therefore, for some polynomial $Q$, on gets
\begin{equation}\label{eq: pre Pila Wilkie bound}
\#X(\QQ,T)\le Z_{\lfloor \log T\rfloor}(f) K Q(\lfloor \log T\rfloor).
\end{equation}
Moreover, in this situation when $f$ has a
B\'ezout bound polynomial in $d$, one obtains the following
improvement of the general bound (\ref{eq: Pila Wilkie bound})
\begin{equation}\label{eq: Poly log Pila Wilkie bound}
\exists \beta, \exists\alpha>0,  \forall T\ge 1,  \
\#X(\QQ,T)\le \beta \log^\alpha T.
\end{equation}
Recently several results appeared, establishing  in different cases
bounds for $\#X(\QQ,T)$
as in (\ref{eq: Poly log Pila Wilkie bound}), some of them proving the existence of convenient
parame\-tri\-zations for certain families of sets $X$ with respect to $\log$-bounds as in
(\ref{eq: Poly log Pila Wilkie bound}), the others proving polynomial B\'ezout bounds in some particular cases (see, among these results,  \cite{Be11}, \cite{Be14}, \cite{BinNovarXiv2},
\cite{BinarXiv}, \cite{BinNovarXiv1}, \cite{BoJo15b},  \cite{BoJo15} , \cite{ClPiWi2016},  \cite{CoMi2016}, \cite{JoMiTh11}, \cite{JoTh12} and \cite{Ma11}).
In the same spirit we give hereafter direct Diophantine applications of the polynomial
B\'ezout bounds obtained in previous sections of the paper.

Let $f$ be an analytic function converging on $D_1$ (on $D_8$
for condition $3$ of Theorem \ref{thm: bound for rational points of bounded height }) and let us denote by $\Gamma_f$ its graph
over $D_{\frac{1}{4}}$.
As a consequence of Theorems
\ref{thm: Bound through the transcendency index sequence},
\ref{thm: Bezout  Bound through polynomial hypertranscendence},
\ref{thm: Bezout Bound for lacunary series},
\ref{thm:ODE}
and Corollary
\ref{cor: Bezout at random},
one has
%
\bt\label{thm: bound for rational points of bounded height }
Assume that one of the following conditions is satisfied
\begin{enumerate}
\item  The Taylor coefficient of $f$ at the origin are rational and the growth conditions
(\ref{eq:growth.cond}) are satisfied,
\item $f$ has rational  Taylor coefficient  at the origin satisfying the growth condition
(\ref{eq:growth.cond}) and $\eta_d$ is polynomially bounded,
\item $f(z)=\sum_{k\ge 1} a_k z^{n_k}$, the lacunarity condition (\ref{eq: lacunarity bounded from above}) is satisfied and for some $p>0$, for any $k\ge 1$, $\vert a_k \vert \ge e^{-n_k^p} $,
\item $f$ is a solution of an algebraic differential equation with rational coefficients and initial conditions,
\item
\label{rem:Log bound for random series}
$f$ is a random series, in the sense of Section \ref{Sec: Bautin determinant for random series}.
\end{enumerate}

Then there exist $\alpha, \beta >0$ such that
$$\# \Gamma_f(\QQ,T) \le \beta \log^\alpha T. $$
\et

\br\label{rmk: Pila Wilkie is asymptotically sharp}
Not only for functions definable in some o-minimal structures, but also
for analytic functions, the
asymptotics of (\ref{eq: Pila Wilkie bound}) is sharp, since, for instance by \cite[Example 7.5]{Pi2004}, \cite{Su2002} or \cite{Su2006}, there exist functions analytic on a neighbourhood of a compact interval having asymptotically as many as possible rational points of height at most $T$ in their graph  with respect to (\ref{eq: Pila Wilkie bound}). For instance, for any $\epsilon\in ]0,1[$, more than $\frac{1}{2}e^{2\log^{1-\epsilon}T}$ points, for an infinite sequence of heights $T$.
In consequence one cannot expect polynomial B\'ezout bounds in all degree $d$ for these analytic functions, since by (\ref{eq: pre Pila Wilkie bound}) one has
$$
\frac{1}{2}e^{2\lfloor \log T\rfloor^{1-\epsilon}}
\le
\frac{1}{2}e^{2\log^{1-\epsilon}T}
\le
Z_{\lfloor \log T\rfloor}(f)  C_{\Gamma_f, {\lfloor \log T\rfloor}} T^{\tau_{\lfloor \log T\rfloor}}.$$
And thus for any $\zeta\in ]0,1[$, there exists a sequence of degrees $d$ going to infinity and such that
\begin{equation}\label{eq: lot of zeroes}
\cal Z_d(f)\ge Z_d(f)\ge  e^{d^\zeta}.
\end{equation}
As a consequence of (\ref{eq: lot of zeroes}),
the condition that $\Delta_d\ge e^{-U(d)}$ for some positive polynomial $U$, may not be satisfied for particular analytic functions $f$. Indeed,
when $\Delta_d\not=0$, the transcendence index $\nu_d$ is polynomially bounded in $d$ and  in case $\Delta_d\ge e^{-U(d)}$,
by Proposition
\ref{pr: Bezout and index of transcendency}
$$\cal Z_d(f)\le 5\nu_{2d}\log(4+(\nu_{2d}+1)\frac{e^{2(d+1)^3}}{\Delta})$$
$$\le
5\nu_{2d}[ \log(5)+\log(\nu_{2d}+1)+2(d+1)^3+U(d)].$$

\er

\br
The condition $\vert a_k \vert \ge e^{-n_k^p}$ of Theorem
\ref{thm: bound for rational points of bounded height } $(4)$ allows
in particular order $0$ for the lacunary series $f(z)=\sum_{k\ge 1} a_k z^{n_k}$ when $f$ is an entire function, since the order of $f$
is given by  $ \limsup_{n\to \infty} -\frac{n_k\log n_k}{\log \vert a_k \vert}$ (see \cite[Theorem 14.1.1]{Hille}), contrariwise
to \cite[Theorem 1.1]{BoJo15b} and \cite[Section 7]{Com.Pol2007} where order $0$ is not allowed. Furthermore the conditions of Theorem \ref{thm: bound for rational points of bounded height } allow to consider analytic functions that are not entire.
\er

\br
Statement \eqref{rem:Log bound for random series} of Theorem  \ref{thm: bound for rational points of bounded height } can be seen as a consequence of
\cite[Theorem 2.7]{BoJo15}, that in fact shows that in case $f(0)$ has a convenient
transcendence measure, then the set $\#\Gamma_f(\QQ,T)$ satisfies the conclusion of  Theorem  \ref{thm: bound for rational points of bounded height }, and on the other hand the set of real numbers having this convenient
transcendence measure is a full set of $\RR$.
\er

%
\br
Using the estimates of \cite[Theorem 2.20]{CoMi2016} (see also
\cite[proof of Theorem 1.5]{Pi2006}) and Corollary \ref{cor: Bezout at random}, on deduces that for random series,
with probability one, the exponent $\alpha$ in the bound of Theorem \ref{thm: bound for rational points of bounded height }
may be chosen as $8$.
\er

%
\br
The conditions $1$ to $4$  on $f$ in  Theorem \ref{thm: bound for rational points of bounded height }
are natural with the aim of showing that there are few rational points in $\Gamma_f$,
since combinations of finer conditions
are considered in order to obtain more remarkable Diophantine properties for $\Gamma_f$.
For instance, in Siegel-Shidlovskii's theorem a combination of conditions comparable to conditions
$1$ and $4$ of Theorem \ref{thm: bound for rational points of bounded height },
among others, imply that there is at most one rational point in $\Gamma_f$.
More accurately, let  $f$ be a $E$-function, that is a Taylor series
$\sum_{k=0}^\infty a_kz^k$ with, for simplicity,  rational coefficients $a_k$ satisfying for any $\varepsilon>0$

$$\vert a_k \vert  =
 O (    k^{k(\varepsilon -1)}   )
\ \hbox{ and } \  \vert q_k \vert =O(k^{\varepsilon k}), $$
where $q_k$ is a common denominator for $a_0, a_1, 2!a_2, \ldots, k!a_k$. Assuming moreover that $f$ is solution of a linear differential equation of order $n$ with coefficients in
$\QQ[z]$,
such that $f, f', \ldots, f^{(n-1)}$ are algebraically independent over $\CC(z)$, then for any algebraic
number $z_0\not=0$, the numbers $f(z_0), \ldots, f^{(n-1)}(z_0)$ are algebraically independent over $\QQ$, and in particular $f(z_0)$ is transcendental
(see for instance  \cite[Theorem 2.1]{Ne2009} or \cite[Theorem 3, p. 123]{Shi1989}).
\er





\bibliographystyle{acm}


\bibliography{Z}



\end{document}